\documentclass[11pt]{article}
\usepackage{amsmath,amssymb,graphics,amsthm}

\newtheorem{lemma}{Lemma}
\newtheorem{remark}{Remark}
\newtheorem{proposition}{Proposition}
\newtheorem{question}{Question}

\newtheorem{theorem}{Theorem}

\begin{document}
\title{{\sc The link surgery of $S^2\times S^2$ and Scharlemann's manifolds}}
\author{Motoo Tange}
\date{}
\maketitle
\abstract
Fintushel-Stern's knot surgery gave many pairs of exotic manifolds, which are
homeomorphic but non-diffeomorphic.
We show that if an elliptic fibration has two parallel, oppositely oriented vanishing circles
(for example $S^2\times S^2$ or Matsumoto's $S^4$), then
the knot surgery gives rise to standard manifolds.
The diffeomorphism can give an alternative proof that Scharlemann's manifold is standard (originally by Akbulut \cite{[Ak1]}).
\section{Introduction.}
\subsection{Knot surgery.}
\hspace{.5cm}Let $X$ be a 4-manifods containing the cusp neighborhood $C$,
which is the well-known elliptic fibration over $D^2$
with one cusp singularity.
In \cite{[FS]} R. Fintushel and R. Stern constructed exotic structures by performing 
the knot surgery of a general fiber of $X$ near the cusp fiber.
The knot surgery on $X$ that containing $C$ is defined as follows.
Let $K$ be a knot in $S^3$.
For a general fiber $T$ of $C$ the surgery
$$X_{K}:=[X-\nu(T)]\cup_{\varphi_0} [(S^3-\nu(K))\times S^1]$$
is called {\it (Fintushel-Stern's) knot surgery}.
Here $\nu(\cdot)$ stands for the interior of the tubular neighborhood.
The gluing map 
$$\varphi_0: \partial\nu(K)\times S^1\to \nu(T^2)=T^2\times\partial D^2$$
 satisfies the following.
$$\text{the meridian of $K$}\times \{\text{pt}\},\ \{\text{pt}\}\times S^1\to\alpha,\ \beta ,$$
$$\text{the longitude of $K$}\times \{\text{pt}\}\to \{\text{pt}\}\times \partial D^2$$
where $\alpha,\beta$ are generators of $H_1(T)$.

We can easily check that $X_{K}$ is homeomorphic to $X$ from Freedman's celebrated result if $X$ is simply connected and closed.
When is $(X,X_{K})$ an exotic pair?
Fintushel and Stern proved the following formula on the Seiberg-Witten invariant.
\begin{equation}
SW_{X_{K}}=SW_X\cdot\Delta_K,
\label{1}
\end{equation}
where $\Delta_K$ is the Alexander polynomial of $K$.
This formula implies that many knot-surgeries  change the differential structures.
However in the case where $\Delta_K(t)=1$ or $SW_X=0$, it is in general unknown whether the pair is exotic or not.

\begin{figure}[htbp]
\begin{center}
\includegraphics{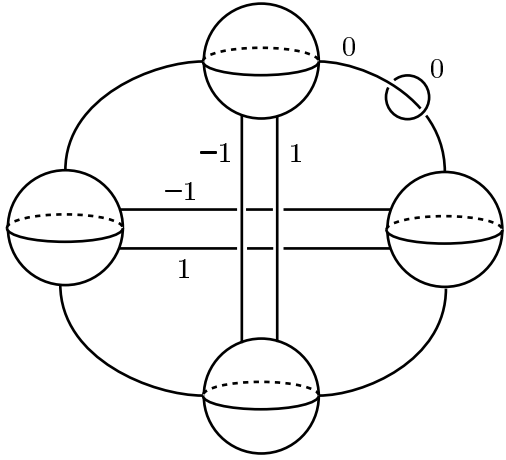}
\caption{The cusp neighborhood and anti-cusp neighborhood in $S^2\times S^2$.}
\label{genfib1}
\end{center}
\end{figure}

On the other hand it is well-known that $S^2\times S^2$ admits achiral Lefschetz fibration containing $C$,
namely $S^2\times S^2$ is diffeomorphic to the double $\overline{C}\cup C$ of $C$.
The diagram is drwan in Figure~\ref{genfib1}.

We denote $\overline{C}\cup C_K$ by $A_K$.
Since $SW_{S^2\times S^2}=0$ holds, SW-invariant cannot distinguish whether $A_K$ is exotic or not. 

In \cite{[Ak2]} S.Akbulut showed that $A_{3_1}$ is diffeomorphic to $S^2\times S^2$.
The diffeomorphism is due to his result \cite{[Ak1]}.
This says the existence of exotic embedding of $C$ into $S^2\times S^2$.
In the article we will show the following.
\begin{theorem}
\label{main}
Let $K$ be any knot.
Then $A_K$ is diffeomorphic to $S^2\times S^2$.
\end{theorem}

\subsection{Link surgery.}
Let $L$ be a link in $S^3$.
We define $A_L$ to be the link surgery of $S^2\times S^2$ along a general fiber near the cusp fiber.
The precise definition of link surgery is in Section~\ref{link}.
Then we give a classification of $A_L$ by applying the same method as the knot case.
\begin{theorem}
\label{links}
Let $L=K_1\cup \cdots\cup K_n$ be any $n$-component link.
Then $A_L$ is diffeomophic to
$$A_L=\begin{cases}\#^{2n-1}S^2\times S^2& \sum_{i\neq j}lk(K_i,K_j)=0\ (2)\ \ \ \text{for any }i\\\#^{2n-1}{\Bbb C}P^2\#^{2n-1}\overline{{\Bbb C}P^2}&\text{otherwise.}\end{cases}$$
\end{theorem}
In Section~\ref{link} we prove the theorem.

\subsection{Scharlemann's manifolds.}
\label{introscha}
We define closed 4-manifolds $B_{K}^{\epsilon}(\gamma)\ (\epsilon=0,1)$ to be the surgery:
$$[S^3_{-1}(K)\times S^1-\nu(\gamma\times \{\text{pt}\})]\cup_{\phi'} S^2\times D^2,$$
where $\gamma$ is a knot in $S^3_{-1}(K)$.
The diffeomorphism type depends only on free homotopy class of the map $S^1\hookrightarrow S^3_{-1}(K)$.

The map $\phi'$ is gluing map $S^2\times \partial D^2\to \partial D^3\times S^1=\partial \nu(\gamma\times \{\text{pt}\})$ where $\phi'(x,t)=(\phi(x,t),t)$.
If $\epsilon=0$, then the gluing map is trivial ($\Leftrightarrow \phi(x,\cdot)=0\in \pi_1(SO(3))$)
and if $\epsilon=1$, then the gluing map is non-trivial ($\Leftrightarrow \phi(x,\cdot)\neq 0$).
In this paper $B_{K}^1(\gamma)$ is called {\it Scharlemann's manifold}.
Any Scharlemann's manifold is homotopy equivalent to $S^3\times S^1\#S^2\times S^2$ (see \cite{[Sc]}).
More strongly $B_{K}^\epsilon(\gamma)$ is homeomorphic to $S^3\times S^1\#S^2\times S^2$ by Freedman's result.
It had been for a long time unknown whether $B_{3_1}^1(\gamma)$ is diffeomorphic to $S^3\times S^1\#S^2\times S^2$ or not.
Akbulut \cite{[Ak1]} showed the following.
\begin{theorem}[\cite{[Ak1]}]
\label{As}
Let $\gamma_0$ be the meridian of $S^3-3_1\subset S^3_{-1}(3_1)$.
$B_{3_1}^{1}(\gamma_0)$ is diffeomorphic to $S^3\times S^1\#S^2\times S^2$.
\end{theorem}

Here we state the following as the third main theorem.
\begin{theorem}
\label{main2}
Suppose that $K$ is any knot.
Let $\gamma_0$ be the meridian of $S^3-K\subset S^3_{-1}(K)$.
Then $B^{1}_{K}(\gamma_0)$ is diffeomorphic to $S^3\times S^1\#S^2\times S^2$.
\end{theorem}

Moreover in Section~\ref{mscha} we will consider the diffeomorphism type of $B_{3_1}^\epsilon(\gamma)$ for other homotopy class $\gamma$.

After I posted the article on arXiv, S. Akbulut also proved the same results in \cite{[Ak4]} in a simple way.
He uses the fact that $S^2\times S^2$ is the double of a cusp neighborhood $C$, however
my proof is available whenever there exist two, parallel, opposite, fishtail singular fiber in an elliptic fibration.

\section*{Acknowledgments.}
The author originally was taught the candidates $A_K$ of exotic $S^2\times S^2$ by Professor M.Akaho (\cite{[Aka]}).
This paper is the negative but complete answer for his question.
I thank him for telling me about attractive 4-dimensional world.

Last year we posed the diffeomorphism for figure-8 surgery of $S^2\times S^2$, however
the proof included a wrong part.
The author deeply apologizes for Differential Topology Seminar at Kyoto University at 2009 because the author explained the wrong proof.
The error was pointed out by Professor S.Akbulut at Michigan state University in 2010, so
the author expresses the gratitude for him.

The author thanks for Professor Masaaki Ue, Dr. Kouichi Yasui and Shohei Yamada for giving me many useful comments by some seminars.
The research is partially supported by JSPS Research Fellowships for Young Scientists (21-1458).

\section{The logarithmic transformation.}
\subsection{Definition and Gompf's result.}
In this section we shall define the logarithmic transformation.
Let $X$ be an oriented 4-manifold and $T\subset X$ a embedded torus with self-intersection $0$.
The surgery 
$$X_{T,p,q,\gamma}=[X-\nu(T)]\cup_\varphi D^2\times T^2,$$
is called {\it logarithmic transformation},
where the gluing map $\varphi:\partial D^2\times T^2\to \partial \nu(T)$ 
satisfies $\varphi(\partial D^2\times \{\text{pt}\})=q(\{\text{pt}\}\times \gamma)+p(\partial D^2\times \{\text{pt}\})$.
It is well-known that the diffeomorphism type of logarithmic transformation depends only on 
the data $(T,p,q,\gamma)$.
The integer $p$ is {\it multiplicity} of the logarithmic transformation, $\gamma$ is the {\it direction} and $q$ is {\it auxiliary multiplicity}.

If $p=1$, then we call $X_{T,1,q,\gamma}$ a $q$-fold Dehn twist of $\partial\nu(T)$ along $T$ parallel to $\gamma$.
\begin{lemma}[Lemma 2.2 in \cite{[G]}]
\label{gompf}
Suppose $N=D^2\times S^1\times S^1$ is embedded in a 4-manifold $X$.
Suppose there is a disk $D\subset X$ intersecting $N$ precisely in $\partial D=\{q\}\times S^1$ for some $q\in\partial D^2\times S^1$,
and that the normal framing of $D$ in $X$ differs from the product framing on $\partial D\subset \partial N$ by $\pm1$ twist.
Then the diffeomorphism type of $X$ does not change if we remove $N$ and reglue it by a $k$-fold Dehn twist of $\partial N$ along
$S^1\times S^1$ parallel to $\gamma=\{q\}\times S^1$.
\end{lemma}

\begin{figure}[htbp]
\begin{center}
\includegraphics{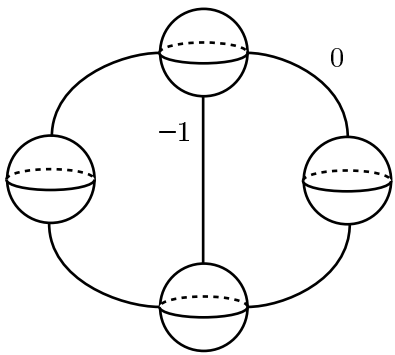}
\caption{Fishtail neighborhood.}
\label{fishtail}
\end{center}
\end{figure}

The manifold $F:=N\cup \nu(D)$ in Lemma~\ref{gompf} is called {\it fishtail neighborhood}, and
the diagram is Figure~\ref{fishtail}.

\subsection{One strand twist.}
\label{onetw}
Let $K_1$ be any knot in $S^3$ and $K_2$ the meridian of $K_1$.
We consider the knot surgery $X_{K_1}$ whose attaching map is $\varphi_0$, where a 4-manifold $X$ contains $C$.
The subset $T_2:=K_2\times S^1\subset X_{K_1}$ is a self-intersection $0$ torus.
Thus the neighborhood $N_2=\nu(K_2)\times S^1$ is the trivial normal bundle over $T_2$.
Since any parallel copy $K_2'\subset \partial N_2$ of $K_2$ by the obvious trivialization of $N_2$
is isotopic to one of vanishing circle of $C_{K_1}$,
there exists a disk $D\subset C_{K_1}$ with $\partial D=K_2'$.
Thus the framing of $\partial D$ coming from the trivialization of $\nu(D)$ differs from the normal framing of the trivialization of $N_2$ by $-1$.
As a result $N\cup \nu(D)$ is the fishtail neighborhood.

From Lemma~\ref{gompf} $n$-fold Dehn twist of $\partial N_2$ along $T_2$ parallel to $K_2$ does not change the diffeomorphism type.
For any integer $n$ we define $\varphi_n$ to be a diffeomorphism satisfying 
\begin{equation}\varphi_n^{-1}(\{\text{pt}\}\times \partial D^2)=[\text{longitude of $K$}]+n[\text{meridian of $K$}],\label{phin}\end{equation}
where other images $\varphi(\text{meridian of $K$})$ and $\varphi(\{\text{pt}\}\times S^1)$ are the same as $\varphi_0$.
We define $X_{K,n}$ to be
$$X_{K,n}:=[X-\nu(T)]\cup_{\varphi_n} [(S^3-\nu(K))\times S^1]$$
The previous paragraph means
$$X_{K_1,n}\cong X_{K_1,0}=X_{K_1}.$$

\begin{figure}[htbp]
\begin{center}
\includegraphics{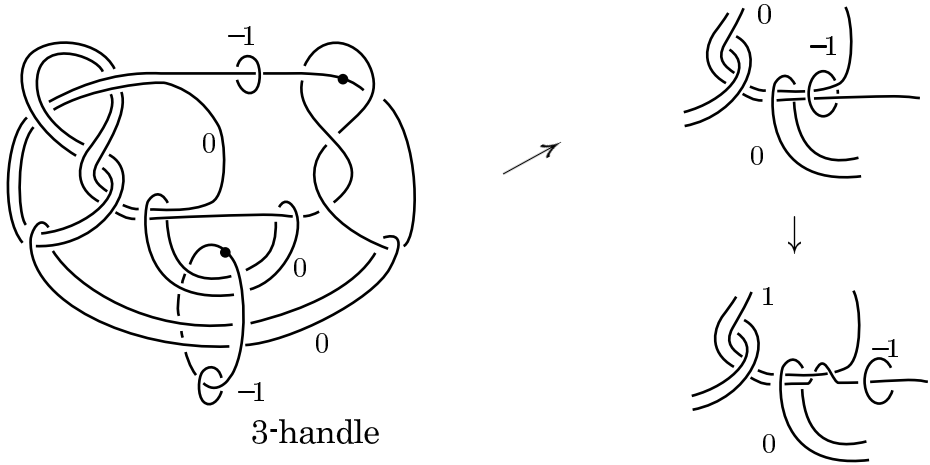}
\caption{The change of framing.}
\label{p6}
\end{center}
\end{figure}

The diffeomorphism can be also understood by handle calculus as Figure~\ref{p6}.
This move is by Akbulut's method in \cite{[Ak1]}.
The left in Figure~\ref{p6} is the $4_1$ surgery of the cusp neighborhood.
Sliding the top $-1$ framed 2-handle over one of two $0$ framed 2-handles below, we get the right-top one in 
Figure~\ref{p6}.
Sliding upper $0$ framed 2-handle over the $-1$ framed 2-handle, we have the right-bottom picture.
This process increases the framing of the knot by $1$.
Iterating the process or the inverse one, we can change the framing to the arbitrary integer.

\section{The knot surgery of $S^2\times S^2$.}
Finding a hidden fishtail neighborhood in $\overline{C}\cup C_K$, we shall prove diffeomorphisms.

{\noindent {\bf Proof of Theorem~\ref{main}.}}
\subsection{Three strand twist.}

\begin{figure}[htbp]
\begin{center}
\includegraphics{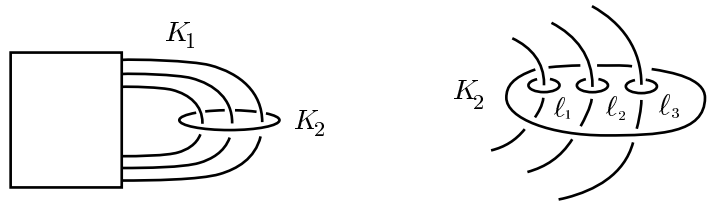}
\caption{$L=K_1\cup K_2$ and $\ell_1,\ell_2,\ell_3$.}
\label{knot3}
\end{center}
\end{figure}

Let $L$ be a two component link as in Figure~\ref{knot3}.
The box is some tangle which presents $K_1$.
We consider the knot surgery $\overline{C}\cup C_{K_1}$ where $\varphi_0$ is the gluing map of the knot surgery.
The torus $T_2=K_2\times S^1\subset [S^3-\nu(K_1)]\times S^1$ has the trivial neighborhood in $A_{K_1}$.
We denote the neighborhood of the torus by $N_2$.

\begin{figure}[htbp]
\begin{center}
\includegraphics{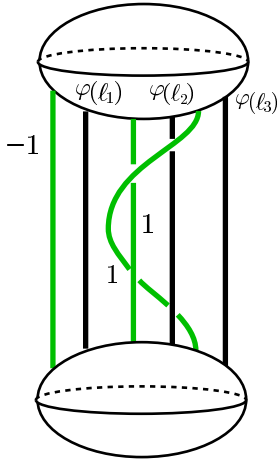}
\caption{An isotopy of $\varphi(\ell_i)$.}
\label{genfib2}
\end{center}
\end{figure}

Our aim here is to construct a fisthtail neighborhood in which $K_2\times S^1$ is a general fiber.
$K_2$ is clearly homologous to the union of meridians $\ell_1,\ell_2$ and $\ell_3$ of $K_1$ via an obvious 3-punctured disk $P$ as in Figure~\ref{knot3}.
Here any $\ell_i$ lies in the boundary of $N_1$ which is the neighborhood of $K_1$.
Each image $\varphi_0(\ell_i)$ is parallel to two vanishing circles in $\overline{C}\cup C_{K_1}$ as in Figure~\ref{genfib2}.
Here the twisted $1$-framed 2-handle is obtained by sliding a trivial 0-framed 2-handle (one of 2-, 3-handle canceling pair) to the $1$-framed 2-handle.

\begin{figure}[htbp]
\begin{center}
\includegraphics{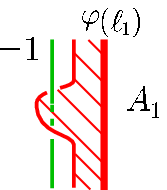}
\caption{$A_1$.}
\label{a1}
\end{center}
\end{figure}

\begin{figure}[htbp]
\begin{center}
\includegraphics{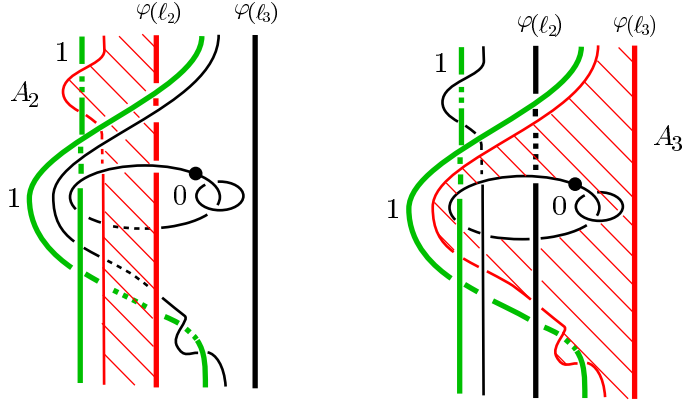}
\caption{Two embedded annuli $A_2$, $A_3$.}
\label{inter}
\end{center}
\end{figure}

Here we will construct three annuli $A_1,A_2$ and $A_3$ that each side of $\partial A_i$ is $\varphi_0(\ell_i)$.
$A_1$ is as in Figure~\ref{a1} and the right side of $\partial A_1$ is $\varphi_0(\ell_1)$.
$A_2$ and $A_3$ are as in the left and right in Figure~\ref{inter} respectively.
$A_3$ runs through the digged 2-handle (the dotted 1-handle) once.
In addition the right sides of $\partial A_2$ and $\partial A_3$ are $\varphi_0(\ell_2)$ and $\varphi_0(\ell_3)$.
Obviously $A_1,A_2$ and $A_3$ are disjoint annuli in $\overline{C}\cup C_{K_1}$.

Another sides of $\partial A_i$ are the boundaries of 2-disks parallel to the cores of the 2-handles in Figure~\ref{genfib2}.
Capping the three 2-disks $C_1$, $C_2$ and $C_3$ to three components of $\partial(P\cup A_1\cup A_2\cup A_3)-K_2$,
we obtain an embedded disk $D:=P\cup A_1\cup A_2\cup A_3\cup C_1\cup C_2\cup C_3$ in $\overline{C}\cup C_{K_1}$
whose boundary is $K_2$.

The framing in $\partial\nu(D)$ coming from the trivialization of $\nu(D)$ differs from the framing of $K_2$ 
coming from the normal bundle of $N_2$ by $-1+1+1=1$.
Therefore $N_2\cup \nu(D)$ is diffeomorphic to $\overline{F}$.
Changing the isotopies of $\varphi_0(\ell_i)$ to two $-1$-framed 2-handles and one $1$-framed 2-handle, we can also embed $F$.


\begin{figure}[htbp]
\begin{center}
\includegraphics{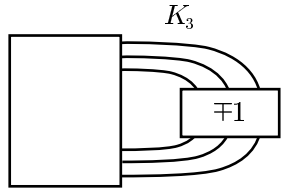}
\caption{$K_3$: a $\pm1$ full twist of $K_1$. The right box is the $\mp1$ full twist.}
\label{dehn}
\end{center}
\end{figure}

Hence $\pm1$-fold Dehn twist along $T_2$ parallel to $K_2$ gives rise to the same manifolds $A_{K_1}$ 
and $\overline{C}\cup C_{K_3,n}$.
The integer $n$ is one of $\mp1$, $\mp9$.
Here $K_3$ is the knot obtained by the $\pm 1$-Dehn surgery along $K_2$ as in Figure~\ref{dehn}.
By one strand twist in Section~\ref{onetw} we have $ A_{K_3}\cong\overline{C}\cup C_{K_3,n}\cong A_{K_1}$.

If we replace three strand twist with odd strand twist, we can construct diffeomorphisms.

\subsection{Ohyama's unknotting operation.}
Y. Ohyama in \cite{[Oh]} has proven that local three strand twist is an unkotting operation of knots.
Therefore for any knot $K$ there exists a finite sequence of local three strand twists: $K=k_0\to k_1\to \cdots \to k_n=\text{unknot}$.
The sequence implies the sequence of diffeomorphisms: 
$$A_{K}=A_{k_0}\cong A_{k_1}\cong \cdots \cong A_{k_n}=S^2\times S^2.$$
\hfill$\Box$

Using the diffeomorphism, we obtain infinitely many embeddings:
$$C\hookrightarrow C\cup\overline{C_{K}}=S^2\times S^2.$$
However whether the embeddings are mutually non-diffeomorphic is unknown.

\begin{remark} 
{\normalfont The diffeomorphism by the three strand twist can be applied to 
any knot surgery of any genus-1 achiral Lefschetz fibration with two oppositely oriented fishtail (or cusp) fiber.
For example $S^4$ admits genus-1 achiral Lefschetz fibration with two 
opposite fishtail fiber (Y. Matsumoto \cite{[M]} and see Figure~8.38 in \cite{[GS]} for the diagram).
Any knot surgery of Matsumoto's Lefschetz fibration on $S^4$ gives rise to standard $S^4$.
}
\end{remark}
\section{Link surgery case.}
\label{link}
Let $L=K_1\cup \cdots,\cup K_n$ be an $n$-component link and $X_i\ (i=1,\cdots,n)$ oriented 4-manifolds which contains the 
cusp neihghborhood $C_i$.
Let $T_i$ be a general fiber of $C_i$.
By the gluing map
$$\varphi_i:\partial\nu(K_i)\times S^1\to \partial \nu(T_i)=T_i\times \partial D^2$$
satisfying
$$\varphi_i(l_i\times \{\text{pt}\})= \{\text{pt}\}\times \partial D^2$$
$$\varphi_i(m_i\times \{\text{pt}\})=\alpha_i,\ \varphi_i(\{\text{pt}\}\times S^1)=\beta_i,$$
where $l_i$ and $m_i$ are the longitude and meridian of $K_i$ and
$\alpha_i$, $\beta_i$ are vanishing circles of $\partial\nu(T_i)$,
we define an operation
$$\coprod_{i=1}^nX_i\to (S^3-\nu(L))\times S^1\cup_{\varphi_i}[X_i-\nu(T_i)].$$
We call the operation {\it link surgery} of $(X_1,\cdots,X_n)$ along $L$ and
denote it by $X(X_1,\cdots,X_n;L)$.
The Seiberg-Witten invariant of $X(X_1,\cdots,X_n:L)$ is computed as follows:
$$SW_{X(X_1,\cdots,X_n:L)}=\Delta_L(t_1,\cdots,t_n)\cdot\prod_{i}^n SW_{E(1)\#_{T=T_i}X_i},$$
where $\Delta_L(t_1,\cdots,t_n)$ is the $n$ variable Alexander polynomial of $L$ and
$E(1)\#_{T=T_i}X_i$ is the fiber sum of the elliptic fibration $E(1)$ and $X_i$ along general fibers $T$ and $T_i$
respectively.
The definition of the fiber sum can be seen in \cite{[FS]}.

\begin{figure}[htbp]
\begin{center}
\includegraphics{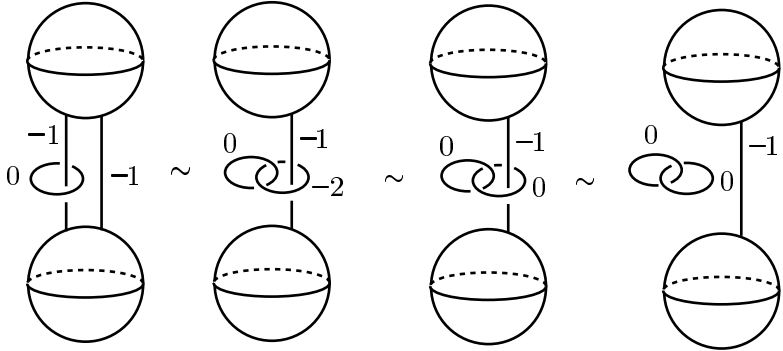}
\caption{$E(1)\#_{T=T_i}S^2\times S^2=E(1)\#^2S^2\times S^2$}
\label{e1s2s2}
\end{center}
\end{figure}

Here we consider the link surgery of $\coprod_{i=1}^nS^2\times S^2$ along any $n$-component link $L$.
We denote the link surgery by $A_L$.
Since the following diffeomorphism
\begin{equation}
\label{sweq}
E(1)\#_{T=T_i}S^2\times S^2\cong E(1)\#^2S^2\times S^2=\#^3{\mathbb C}P^2\#^{11}\overline{{\mathbb C}P^2}
\end{equation}
holds, we have $SW_{A_L}=0$.
The first diffeomorphism is due to Figure~\ref{e1s2s2}.
The leftmost figure is a submanifold of $E(1)\#_{T=T_i}S^2\times S^2$, where 
the handle decomposition of $E(1)$ uses the diagram of Figure 8.10 in \cite{[GS]}.
Sliding handles several times, we find a separated Hopf link (the rightmost figure).
The Hopf link is unlinked from other attaching handles of $E(1)\#_{T=T_i}S^2\times S^2$ while
the handle sliding is operated.
The second equality of (\ref{sweq}) holds by some blow ups and downs .
Thus by vanishing theorem of $SW$-invariant, we have $A_{L}=0$.

\begin{figure}[htbp]
\begin{center}
\includegraphics{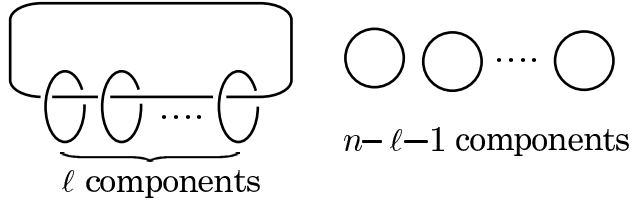}
\caption{The representation $L_{n,\ell}$ of ${\mathcal L}_n$}
\label{repre}
\end{center}
\end{figure}
{\noindent {\bf Proof of Theorem~\ref{links}.}} \
Let $L=K_1\cup K_2\cup \cdots \cup K_n$ be any $n$-component link.
The set $\tilde{{\mathcal L}}_n$ of all $n$-component links up to local three strand twist consists of $2^{n-1}$ classees
due to Nakanishi and Ohyama's results \cite{[Oh],[Na]}.
Forgetting the ordering of the components of any link in $\tilde{{\mathcal L}}_n$ we get
a set ${\mathcal L}_n$.
The set ${\mathcal L}_n$ has $n$ classes.
The class is represented by $L_{n,\ell}\ (\ell=0,1,\cdots,n-1)$ as in Figure~\ref{repre}.

By applying three strand twist method to link surgery $A_L$ we have only to consider diffeomorphism type of $A_{L_{n,\ell}}$ for some $\ell$.

\begin{figure}[htbp]
\begin{center}
\includegraphics{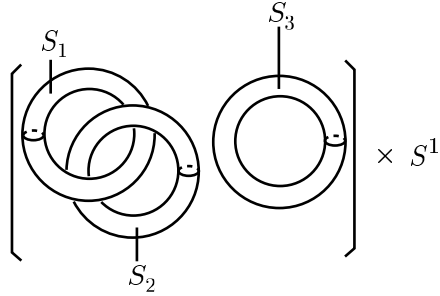}
\caption{Torus neighborhoods $S_i\times S^1$ in $S^3\times S^1$.}
\label{sl}
\end{center}
\end{figure}

\begin{figure}[htbp]
\begin{center}
\includegraphics{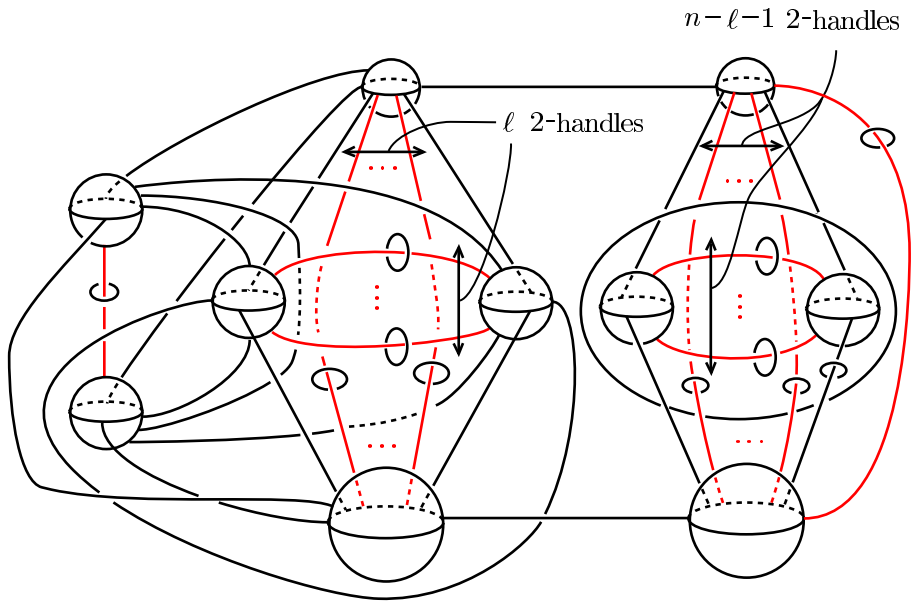}
\caption{$\coprod_{i=1}[\overline{C}\cup C-\nu(T)]\cup[(S^3-\nu(L_{n,\ell}))\times S^1]$.}
\label{pr1}
\end{center}
\end{figure}

Now suppose that $\ell\ge 1$.
The $\ell$ and $n-\ell-1$ components in Figure~\ref{repre} are parallel each other.
Then we embed $T^2\times D^2$ in $S^3\times S^1$ as $S_i\times S^1$ $(i=1,2,3)$ where $S_1$, $S_2$, and $S_3$ are embedded solid tori in $S^3$ as in Figure~\ref{sl}.
Let $C(m)$ be elliptic fibration over $D^2$ with exactly $m$ parallel cusp fibers as singular fibers.
$A_{L_{n,\ell}}$ is obtained by attaching $\overline{C(1)}\cup C(1)-\nu(T_1)$ to $\partial(\nu(S_1)\times S^1)$,
$\overline{C(\ell)}\cup C(\ell)-\nu(T_2)$ to $\partial (\nu(S_2)\times S^1)$ and 
$\overline{C(n-\ell-1)}\cup C(n-\ell-1)-\nu(T_3)$ to $\partial(\nu(S_3)\times S^1)$.
The tori $T_1$, $T_2$ and $T_3$ are general fibers in the cusp neighborhoods.
Thus the diagram is Figure~\ref{pr1}.

\begin{figure}[htbp]
\begin{center}
\includegraphics{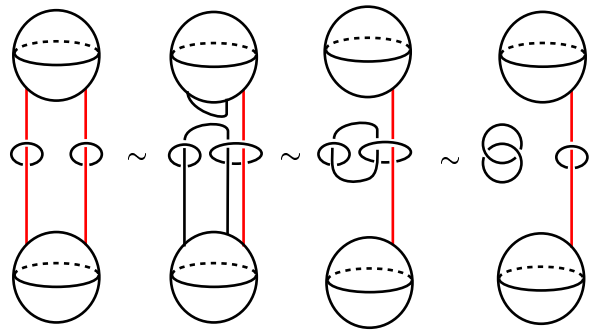}
\caption{To make $S^2\times S^2$-component from two parallel $-1$ framed 2-handles.}
\label{pr2}
\end{center}
\end{figure}

\begin{figure}[htbp]
\begin{center}
\includegraphics{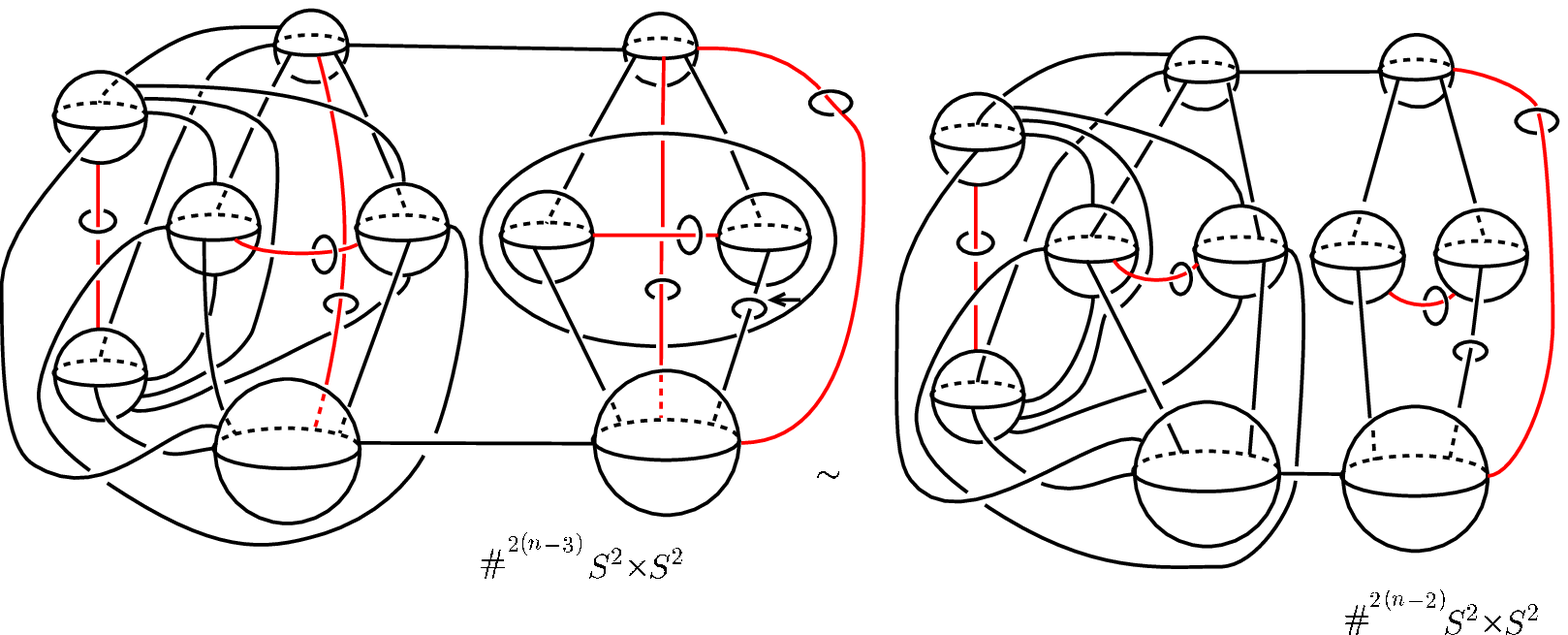}
\caption{To come out two $S^2\times S^2$ components.}
\label{pr3}
\end{center}
\end{figure}

\begin{figure}[htbp]
\begin{center}
\includegraphics{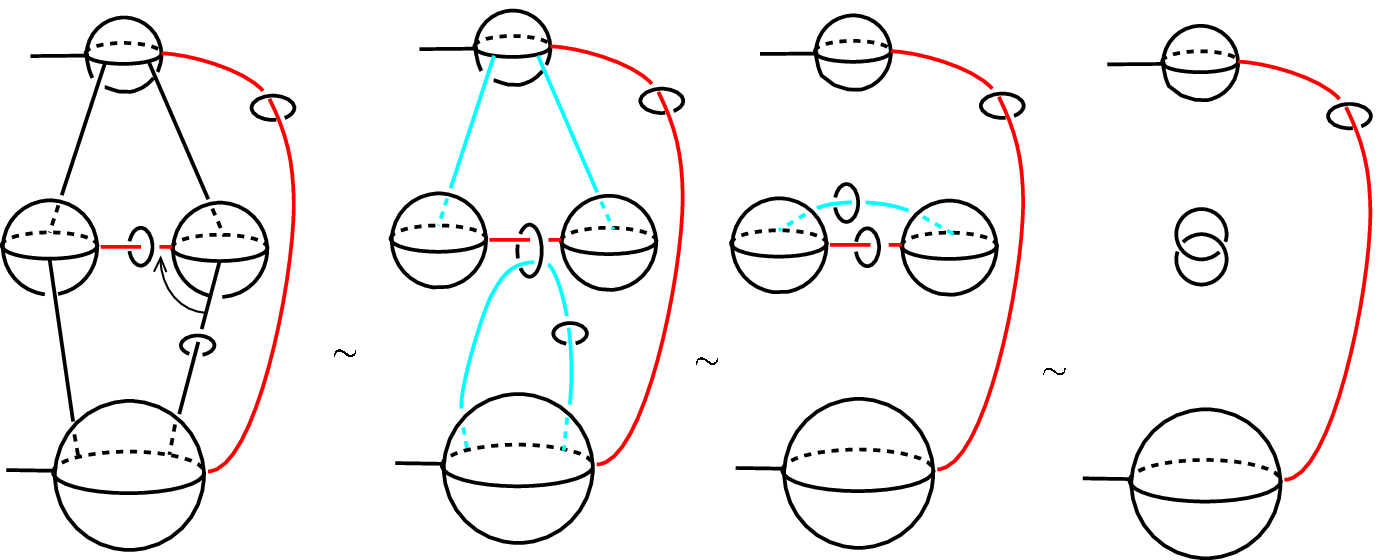}
\caption{The right part of Figure~\ref{pr3}.}
\label{pr4}
\end{center}
\end{figure}

\begin{figure}[htbp]
\begin{center}
\includegraphics{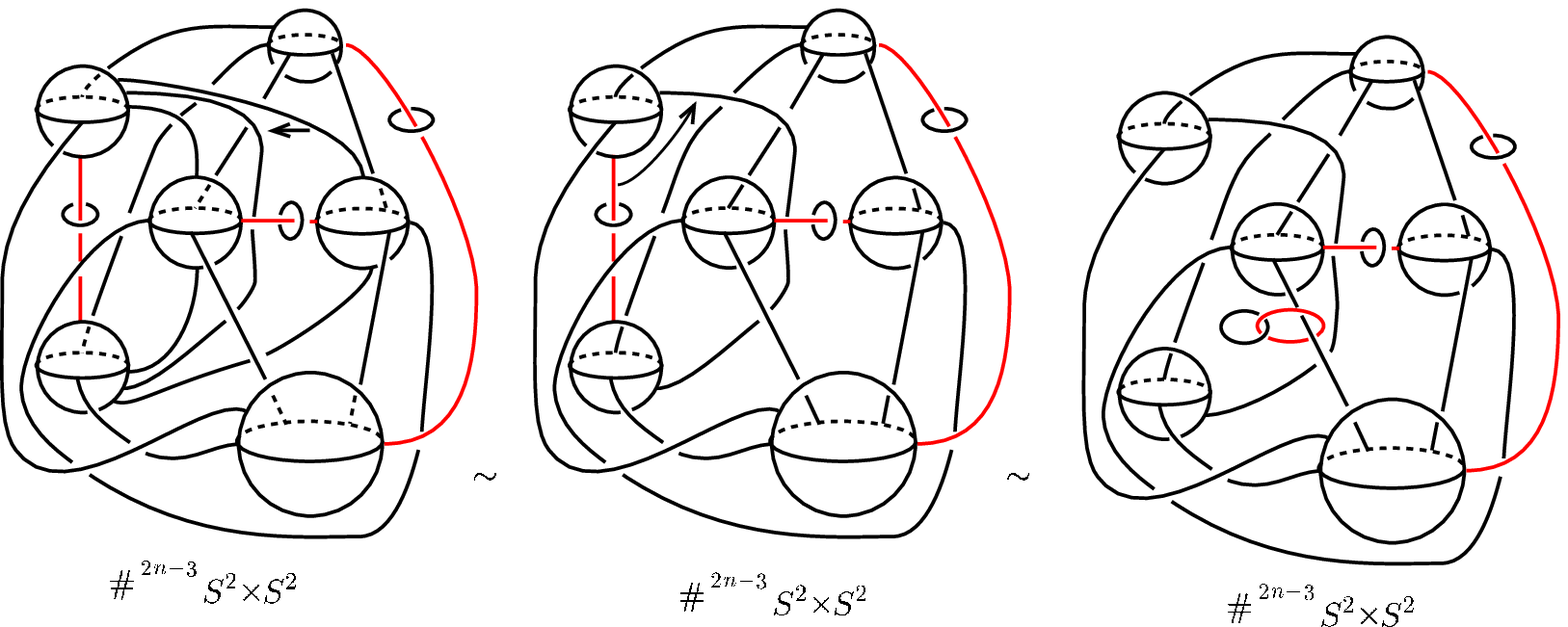}
\caption{Several handle slidings.}
\label{pr5}
\end{center}
\end{figure}

\begin{figure}[htbp]
\begin{center}
\includegraphics{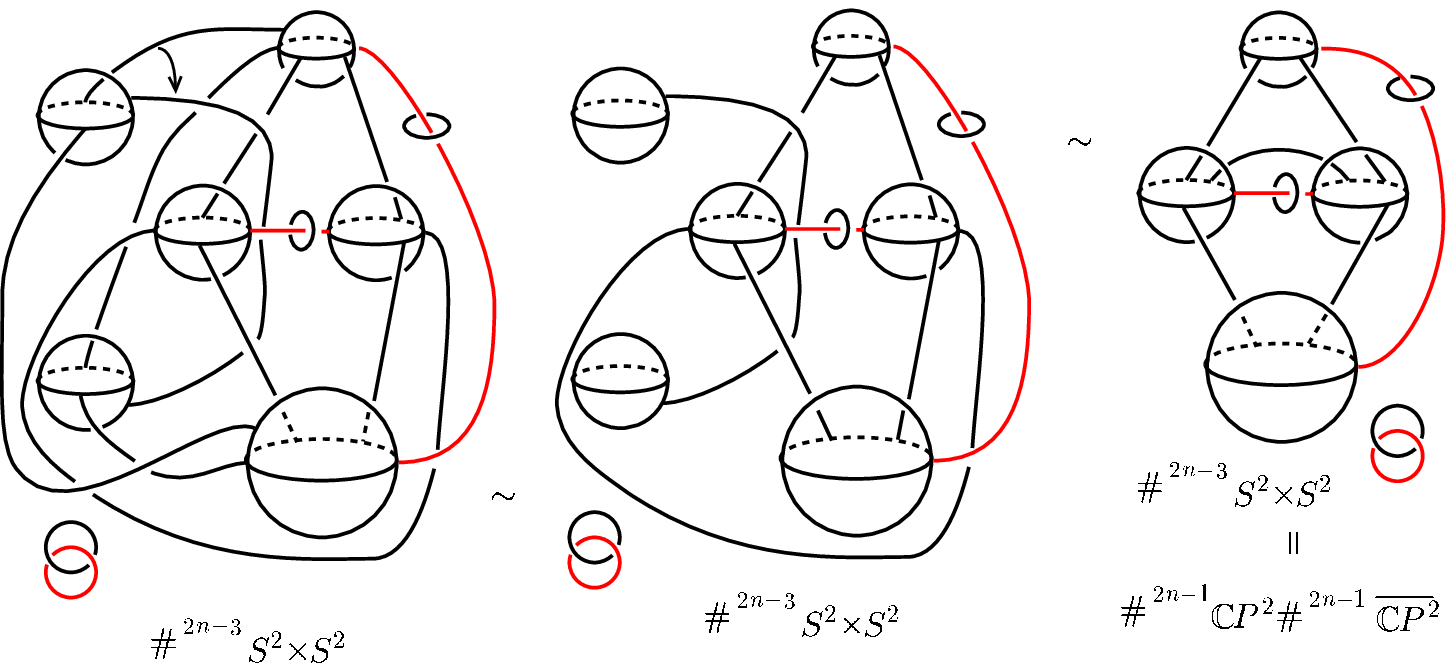}
\caption{}
\label{pr6}
\end{center}
\end{figure}

In this section from here on in the diagram red, black, blue components stand for 
$-1$-, $0$-, $1$-framed 2-handles respectively.
Here parallel $-1$ framed 2-handles with $0$ framed 2-handles in Figure~\ref{pr2} gives rise to $S^2\times S^2$ connected sum components in such a way as Figure~\ref{pr2}.
Moreover we can take out two $S^2\times S^2$ connected components by some handle slidings (Figure~\ref{pr3}).
The right part of Figure~\ref{pr3} are moved as in Figure~\ref{pr4}.
Sliding several times as in Figure~\ref{pr5}-\ref{pr6} we get the right diagram.
The Hopf-linked handles of the right in Figure~\ref{pr6} is a $\Bbb{C}P^2\#\overline{\Bbb{C}P^2}$-component.
The rest of diagram in the left of Figure~\ref{pr6} is another $\Bbb{C}P^2\#\overline{\Bbb{C}P^2}$.
Thus we have $\#^{2n-3}S^2\times S^2\#^2\Bbb{C}P^2\#^2\overline{\Bbb{C}P^2}=\#^{2n-1}\Bbb{C}P^2\#^{2n-1}\overline{\Bbb{C}P^2}$.

\begin{figure}[htbp]
\begin{center}
\includegraphics{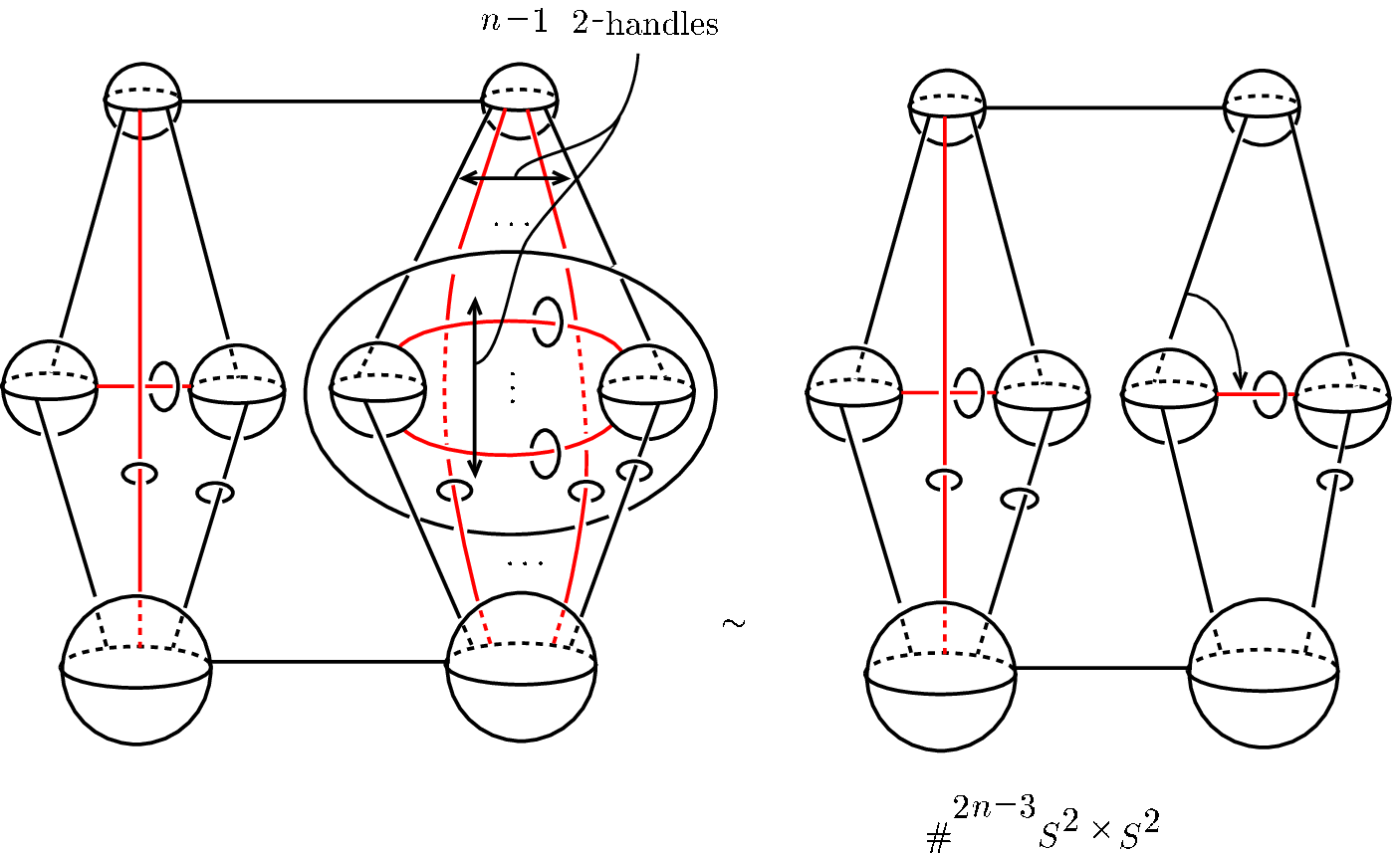}
\caption{The $\ell=0$ case.}
\label{pr7}
\end{center}
\end{figure}

\begin{figure}[htbp]
\begin{center}
\includegraphics{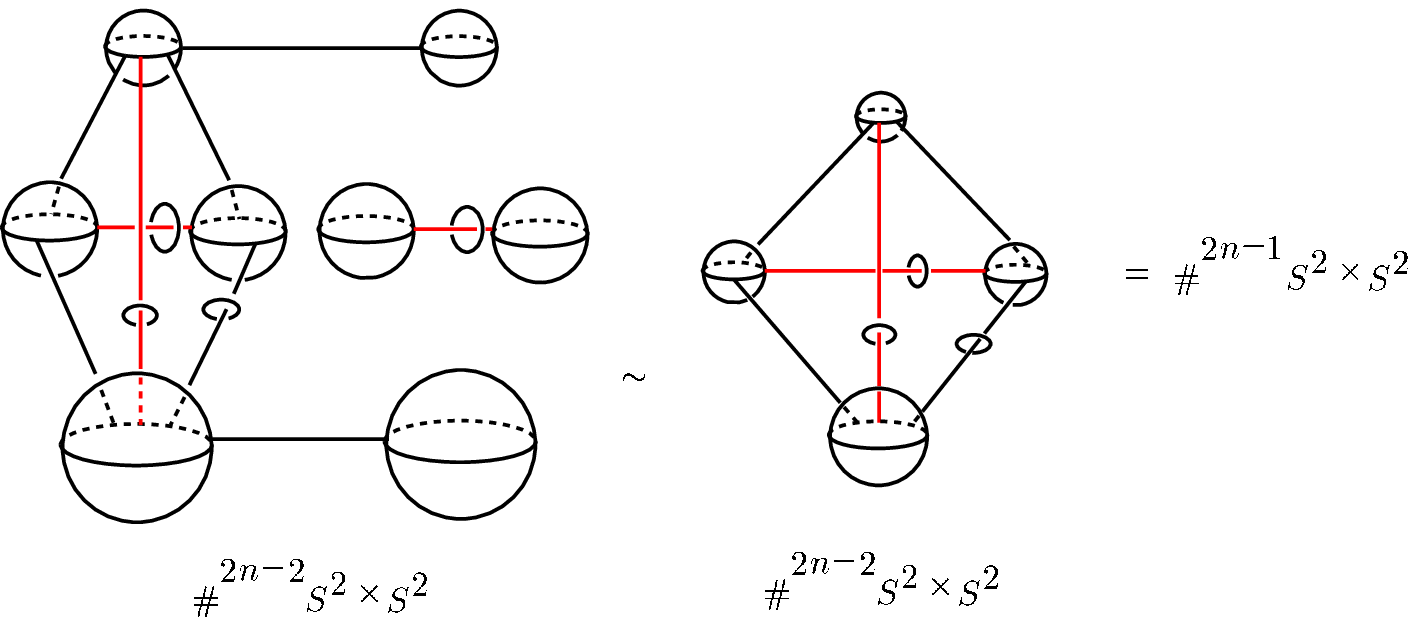}
\caption{Proof of Theorem~\ref{links}.}
\label{pr8}
\end{center}
\end{figure}

Suppose that $\ell=0$.
The diagram of $A_{L_{n,0}}$ is the left of Figure~\ref{pr7}.
Applying Figure~\ref{pr2}, we get the right of Figure~\ref{pr7}.
Sliding handle as indicated in Figure~\ref{pr7} gives rise to the left of Figure~\ref{pr8}.
This sliding is the same as Figure~\ref{pr4}.
Canceling handles in Figure~\ref{pr8}, we obtain $A_{L_{n,0}}\cong \#^{2n-1}S^2\times S^2$.\hfill$\Box$

\section{Scharlemann's manifold.}
From here we focus on Scharlemann's manifolds which are defined in Section~\ref{introscha}.
\subsection{Scharlemann's manifolds along the meridians.}
\begin{figure}[htbp]
\begin{center}
\includegraphics{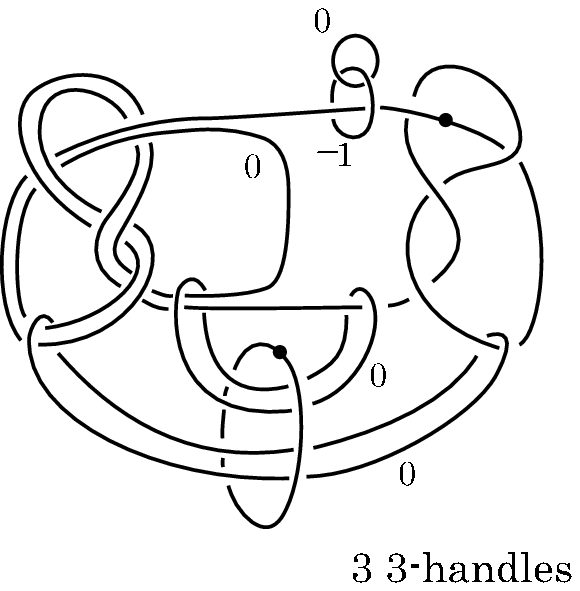}
\caption{$\overline{F}\cup[F-\nu(T)]\cup_{\varphi_{0}} [(S^3-\nu(K))\times S^1]$.}
\label{p7}
\end{center}
\end{figure}

In this subsection we consider Scharlemann's manifolds by the meridian $\gamma_0$.
We remark the following.
\begin{remark}
\begin{normalfont}
Let $\gamma_0$ be the meridian circle in $S^3_{-1}(K)$.
Scharlemann's manifolds $B^0_{K}(\gamma_0)$ is 
always diffeomorphic to $S^3\times S^1\#{\Bbb C}P^2\#\overline{{\Bbb C}P^2}$ by handle calculus.
\end{normalfont}
\end{remark}
In the case of $\epsilon=1$, we note the relationship between $B_K^1(\gamma_0)$ and 
the knot surgery of the fishtail neighborhood.
\begin{lemma}
\label{lem10}
$B_K^1(\gamma_0)$ is diffeomorphic to $\overline{F}\cup F_{K}$. 
\end{lemma}
The gluing map $\varphi_0$ satisfies (\ref{phin}) and $\varphi_0(\text{meridian of $K$})=$ the vanishing cycle.
The resulting manifold 
$$\overline{F}\cup F_K=\overline{F}\cup [F-\nu(T)]\cup_{\varphi_{0}} [(S^3-\nu(K))\times S^1]$$
is Figure~\ref{p7},
where the picture is the case of $K=4_1$.

{\noindent{\bf Proof of Lemma~\ref{lem10}.}
The surgery along $\gamma_0\times \{\text{pt}\}$ in $S^3_{-1}(K)\times S^1$ is Figure~\ref{p8}.
Hence we get the following diffeomorphisms.
\begin{eqnarray}
B_K^1(\gamma_0)&=&[S^3_{-1}(K)\times S^1-\nu(\gamma_0)]\cup_{\phi'} S^2\times D^2=\overline{F}\cup(F-\nu(T))\cup_{\varphi_{-1}}[S^3-\nu(K)]\times S^1\nonumber\\
&\cong &\overline{F}\cup(F-\nu(K))\cup_{\varphi_{0}}S^2\times D^2=\overline{F}\cup F_{K}\nonumber
\end{eqnarray}
\hfill$\Box$\\
\begin{figure}[htbp]
\begin{center}
\includegraphics{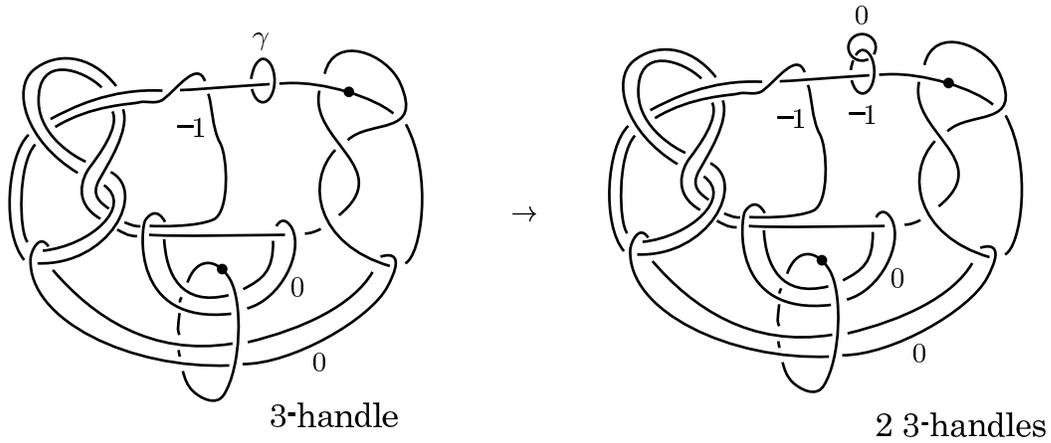}
\caption{The non-trivial surgery along $\gamma$ in $S^3_{-1}(K)\times S^1$.}
\label{p8}
\end{center}
\end{figure}

{\noindent{\bf Proof of Theorem~\ref{main2}.}}
Applying the same diffeomorphism argument as the proof of Theorem~\ref{main} and
Ohyama's unknotting operation imply 
$$\overline{F}\cup F_{K_1}=\overline{F}\cup F\cong S^3\times S^1\#S^2\times S^2.$$
The last diffeomorphism is Figure~\ref{dfishtail}.
\hfill$\Box$
\begin{figure}[htbp]
\begin{center}
\includegraphics{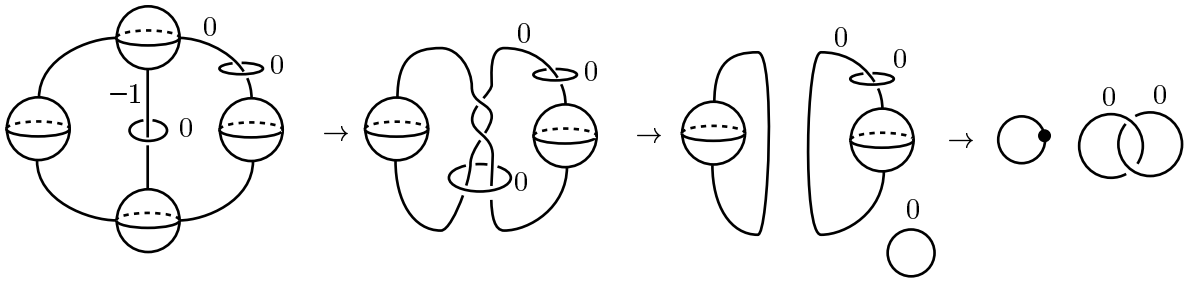}
\caption{The double of $F$.}
\label{dfishtail}
\end{center}
\end{figure}

The manifolds obtained by the same meridian surgery on $S^3_p(K)\times S^1$
are diffeomorphic to $B_K^{\epsilon}(\gamma_0)$ by one strand twist.

\begin{remark}
\begin{normalfont}
$B_K^1(\gamma_0)$ is obtained from $A_K$ as a surgery along an embedded $S^2$.
The neighborhood of the sphere $\Sigma$ is the union of the bottom $0$ framed 2-handle and the 4-handle (the left of Figure~\ref{p3}).
\begin{figure}[htbp]
\begin{center}
\includegraphics{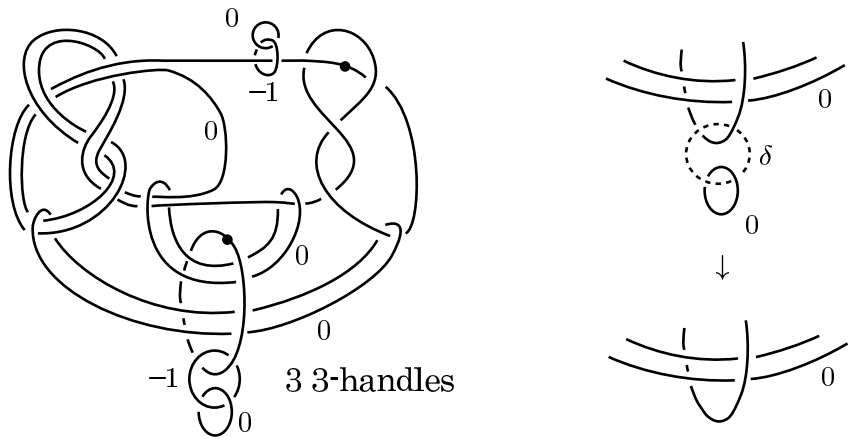}
\caption{The left: $A_K$. The right: surgery $B_K^1=[A_K-\nu(\Sigma)]\cup S^1\times D^3$.}
\label{p3}
\end{center}
\end{figure}
Attaching the 3-handle and 4-handle to the complement gets $B_K^1$ (the right of Figure~\ref{p3}).
The circle $\delta$ in Figure~\ref{p3} is the core circle of $S^1\times D^3$ attached.
\end{normalfont}
\end{remark}

\begin{remark}
\begin{normalfont}
In \cite{[Ak3]} Akbulut got a plug twisting $(W_{1,2},f)$ satisfying $E(1)=N\cup_{\text{id}}W_{1,2}$ and $E(1)_{2,3}=N\cup_{f}W_{1,2}$.
The definitions of plug, $N$ and $W_{1,2}$ are written down in \cite{[Ak3]}.
In the same way as \cite{[Ak3]} we can also show that there exist infinitely many plug twistings $(W_{1,2},f_K)$ of $E(1)$ 
with the same plug $W_{1,2}$.
As a result any plug twisting satisfies $E(1)=M\cup_{\text{id}}W_{1,2}$ and $E(1)_K=M\cup_{f_K}W_{1,2}$.
Infinite variations of Alexander polynomial $\Delta_K(t)$ of knot imply the existence of infinite embeddings $W_{1,2}\hookrightarrow M\cup_{id} W_{1,2}$.
\end{normalfont}
\end{remark}

\subsection{Scharlemann's manifold along a non-meridian circle.}
\label{mscha}
The fundamental group of $S^3_{-1}(3_1)$ is
$$\pi_1(S^3_{-1}(3_1))=\langle x,y|x^5=(xy)^3=(xyx)^2\rangle\cong \tilde{A}_5.$$
The set of free homotopy classes of maps $S^1\to S^3_{-1}(3_1)$
is
\begin{equation}
\label{homotopy}
[S^1,S^3_{-1}(3_1)]=\pi_1(S^3_{-1}(3_1))/\text{conj.}
\end{equation}
and it possesses 9 classes as follows.
\begin{center}
\begin{tabular}{|c|c|c|c|c|c|c|c|c|c|}
\hline
\text{Classes} & $[e]$ & $[x^5]$ & $[xyx]$ & $[x]$ & $[x^2]$ & $[x^3]$ & $[x^4]$ & $[xy]$ & $[(xy)^2]$\\
\hline
\text{Order}&1&2&4& 10 & 5 & 10& 5& 6 & 3 \\
\hline 
\end{tabular}
\end{center}
Each of the classes is a normal generator of the fundamental group except $[e],[x^5]$.
Let $K$ be a knot in $S^3$.
For $z\in \pi_1(S^3_{-1}(3_1))$ we define $B_{3_1}^\epsilon(z)$ to be the surgery of $S^3_{-1}(3_1)\times S^1$
along a knot $K\subset S^3_{-1}(3_1)$ whose free homotopy class $[K]$ coincides with the conjugacy class $[z]$.
Akbulut's result \cite{[Ak1]} is that $B_{3_1}^1(x)$ is the standard $S^2\times S^2\#S^3\times S^1$.
For other conjugacy class we will prove the following.
\begin{proposition}
\label{xy-2}
$B_{3_1}^1(xy)$ is diffeomorphic to $S^3\times S^1\#S^2\times S^2$
and
$B_{3_1}^0(xy)$ is diffeomorphic to $S^3\times S^1\#{\Bbb C}P^2\# \overline{{\Bbb C}P^2}$.
 
\end{proposition}

\begin{remark}{\normalfont
The conjugacy class differs from that of $x$,
for the order of $xy^{-2}$ is $6$ while the order of $x$ is $10$.
Thus the two surgeries are different in general.
}\end{remark}
\begin{figure}[htbp]
\begin{center}
\includegraphics{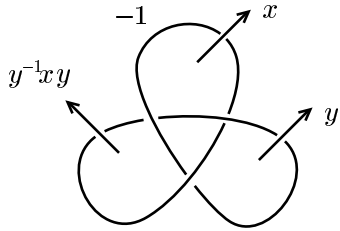}
\caption{The generator $x,y$ of $\pi_1(\Sigma)$.}
\label{tre}
\end{center}
\end{figure}
\begin{remark}{\normalfont
Here taking the diagram in 
Figure~\ref{tre}, 
we fix the framing of $\gamma$.
}\end{remark}

{\noindent{\bf Proof of Proposition~\ref{xy-2}.}}
The meridian $\gamma_0\subset S_1^3(4_1)$ is isotopic to a curve $\gamma_0'\subset S^3_{-1}(3_1)$
as in Figure~\ref{ft}.
The homotopy class $[\gamma_0']$ of the curve is equal to $[xy^{-2}]$ in Equality~(\ref{homotopy}).
In fact $(y^{-1}xy)^2y^{-1}=y^{-1}x^2\sim xy^{-1}x=x^{2}y^{-1}x^{-1}y^{-1}x^2\sim x^5(yxyx)^{-1}=yx\sim xy$ holds.
Hence using one strand twist twice and Theorem~\ref{main2}, we get the following diffeomorphisms
\begin{eqnarray*}
B_{3_1}^1(xy)&=&[S^3_{-1}(3_1)\times S^1-\nu(\gamma_0'\times \{\text{pt}\})]\cup_{\phi'}D^2\times S^2\\
&\cong&[S^3_{1}(4_1)\times S^1-\nu(\gamma_0\times \{\text{pt}\})]\cup_{\phi'}D^2\times S^2\\
&\cong&[S^3_{-1}(4_1)\times S^1-\nu(\gamma_0\times \{\text{pt}\})]\cup_{\phi'}D^2\times S^2\\
&=&B_{4_1}^1(\gamma_0)\cong S^2\times S^2\#S^3\times S^1.
\end{eqnarray*}
\begin{figure}[htbp]
\begin{center}
\includegraphics{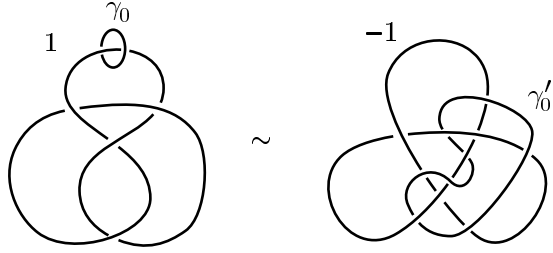}
\caption{An isotopy of $\gamma_0$ via the diffeomorphism $S^3_{1}(4_1)\cong S^3_{-1}(3_1)$.}
\label{ft}
\end{center}
\end{figure}
The case of $B_{3_1}^0(xy)$ is diffeomorphic to $S^3\times S^1\#{\Bbb C}P^2\#\overline{{\Bbb C}P^2}$
in the same way.\hfill$\Box$

Here we will argue several other cases.
\begin{proposition}
$B_{3_1}^1([x^2])$, $B_{3_1}^0([x^3])$ and $B_{3_1}^1([x^4])$ are diffeomorphic to $S^3\times S^1\#{\Bbb C}P^2\#\overline{{\Bbb C}P^2}$.
\end{proposition}
{\bf Proof.}
Here we use notation $\approx$, $\sim$ and $=$ as one strand twist, a homeomorphism, Kirby calculus technique in 3-dimansional.

In the case of $B_{3_1}^1([x^2])$, the last picture in 
Figure~\ref{tr2} 
represents a surgery on $L(5,-1)\times S^1$.
This is diffeomorphic to $S^3\times S^1\#{\Bbb C}P^2\#\overline{{\Bbb C}P^2}$.
\begin{figure}[htbp]
\begin{center}
\includegraphics{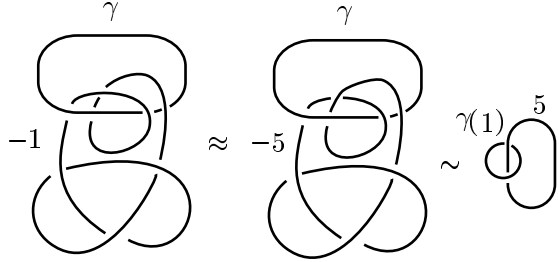}
\caption{The diffeomorphism for $B_{3_1}^1([x^2])$.}
\label{tr2}
\end{center}
\end{figure}

In the case of $B_{3_1}^0([x^3])$ the last picture in 
Figure~\ref{tr3} 
represents a surgery along a knot $\gamma$ in $S^3\times S^1$ with odd framing.
Namely the manifold is $S^3\times S^1\#{\Bbb C}P^2\#\overline{{\Bbb C}P^2}$.
\begin{figure}[htbp]
\begin{center}
\includegraphics{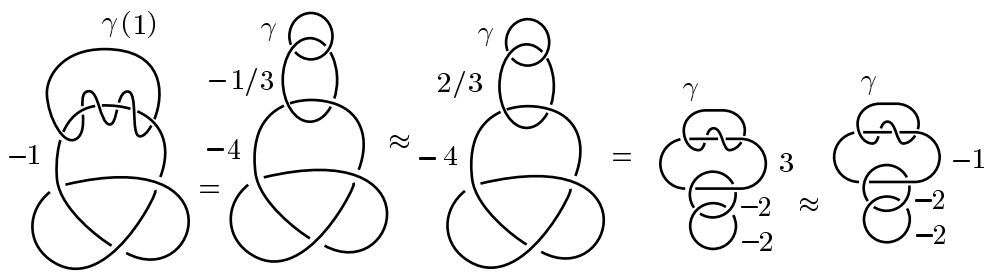}
\caption{The diffeomorphism for $B_{3_1}^0([x^3])$.}
\label{tr3}
\end{center}
\end{figure}

In the case of $B_{3_1}^1([x^4])$, the last picture in 
Figure~\ref{tr4} 
gives $S^3\times S^1\#{\Bbb C}P^2\#\overline{{\Bbb C}P^2}$ in the similar way.
Here $Pr(-2,3,7)$ is the $(-2,3,7)$-pretzel knot.
\begin{figure}[htbp]
\begin{center}
\includegraphics{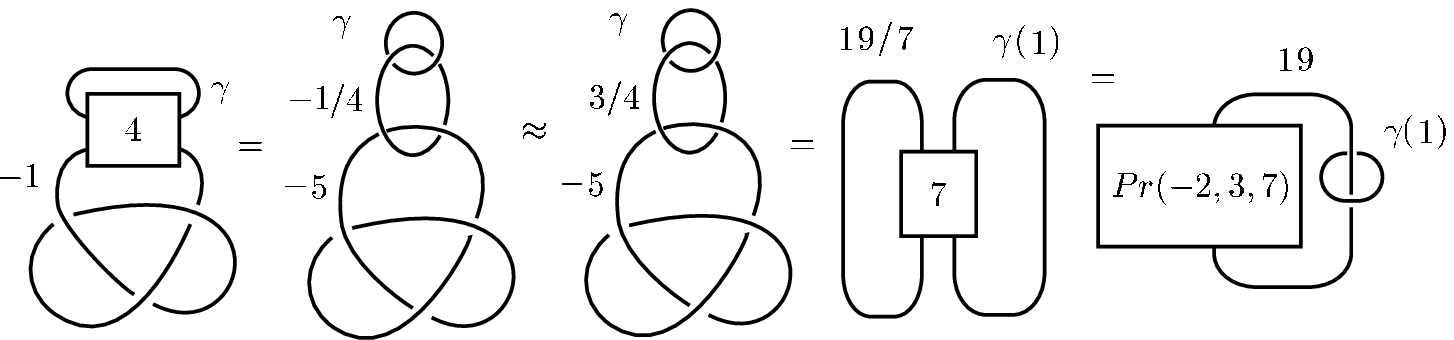}
\caption{The diffeomorphism for $B_{3_1}^1([x^4])$.}
\label{tr4}
\end{center}
\end{figure}
\hfill$\Box$

\begin{proposition}
$B_{3_1}^0(xyx)$ is diffeomorphic to $S^3\times S^1\#{\Bbb C}P^2\#\overline{{\Bbb C}P^2}$.
\end{proposition}
{\bf Proof.}
The homotopy cass of the curve $\gamma$ is
$xyy^{-1}xy\sim x^2y\sim xyx$. then 
from the 
Figure~\ref{tr5}, 
we get $S^3\times S^1\#{\Bbb C}P^2\#\overline{{\Bbb C}P^2}$.
\hfill$\Box$.

\begin{figure}[htbp]
\begin{center}
\includegraphics{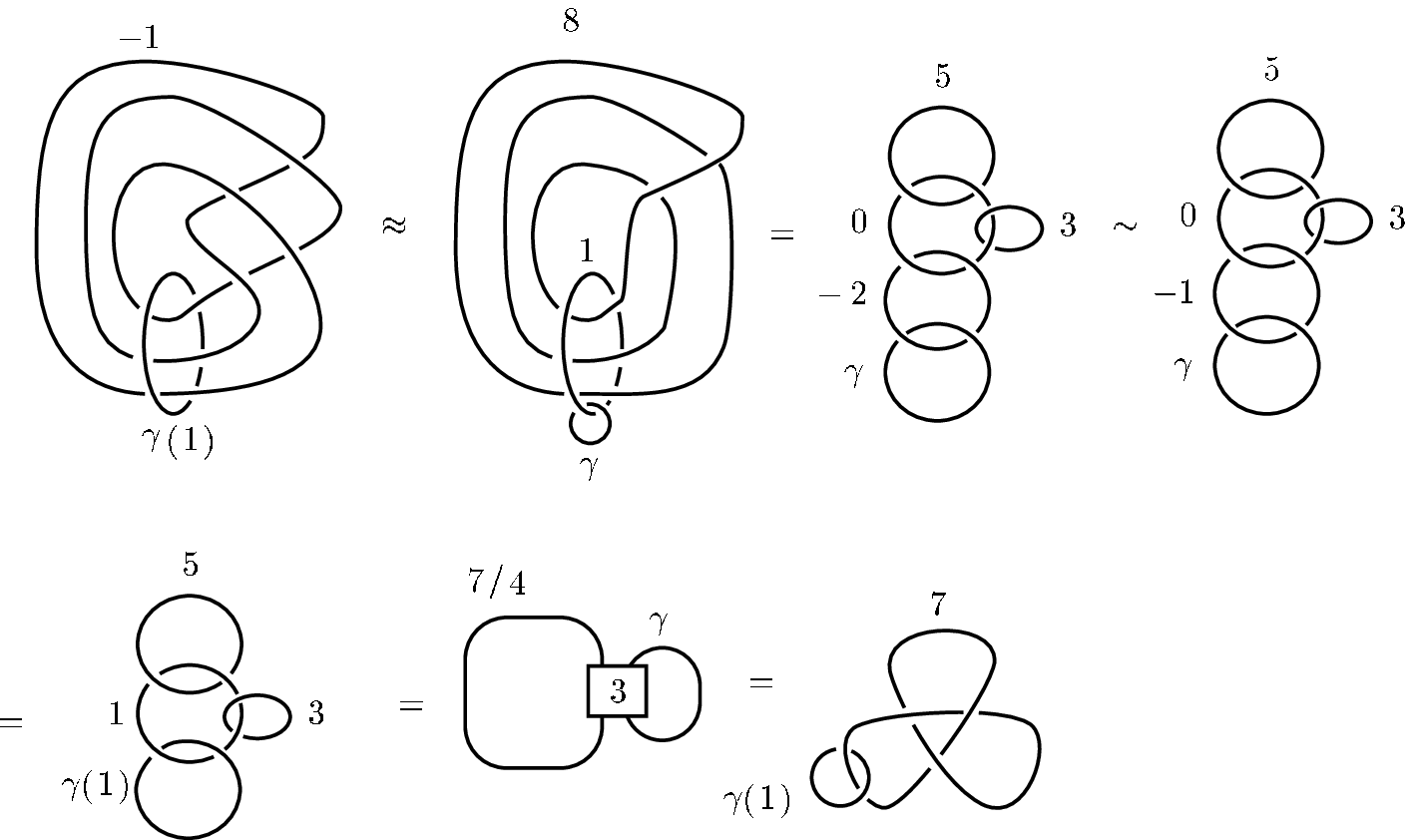}
\caption{The diffeomorphism for $B_{3_1}^0(xyx)$.}
\label{tr5}
\end{center}
\end{figure}

\begin{proposition}
$B_{3_1}^1((xy)^2)$ is diffeomorphic to $S^3\times S^1\#{\Bbb C}P^2\#\overline{{\Bbb C}P^2}$.
\end{proposition}
{\bf Proof.}
The deformation as in 
Figure~\ref{tr6} 
we get $S^3\times S^1\#{\Bbb C}P^2\#\overline{{\Bbb C}P^2}$.
Here $T_{2,7}$ is the $(2,7)$-torus knot
\hfill$\Box$.

\begin{figure}[htbp]
\begin{center}
\includegraphics{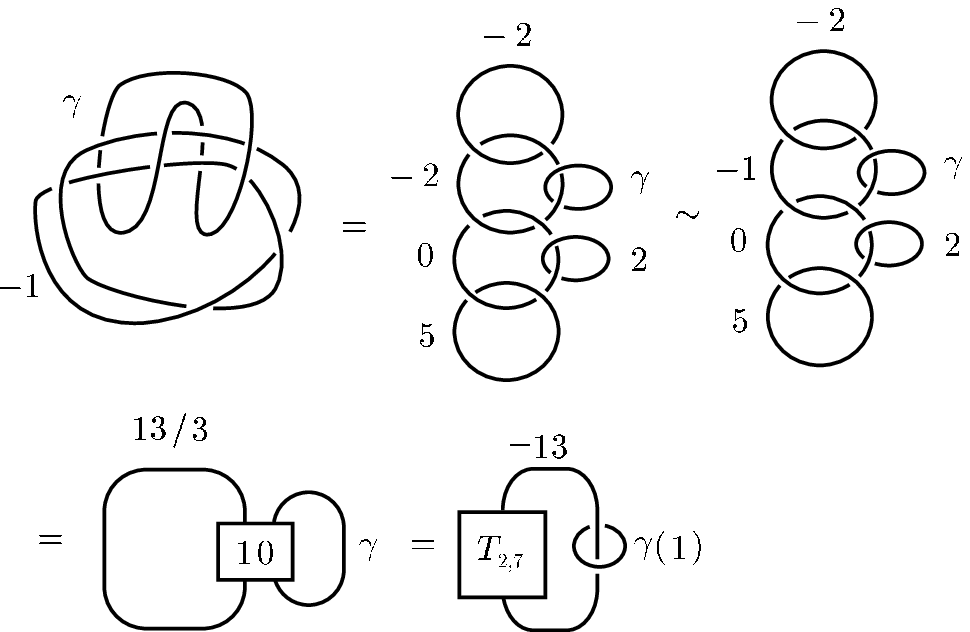}
\caption{The diffeomorphism for $B_{3_1}^1((xy)^2)$.}
\label{tr6}
\end{center}
\end{figure}

In the end we raise the questions.
\begin{question}
Are the following manifolds
$$B_{3_1}^0(x^2),\ B_{3_1}^1(x^3),\ B_{3_1}^0(x^4),\ B_{3_1}^1(xyx),\ B_{3_1}^0((xy)^2)$$
diffeomorphic to $S^3\times S^1\#S^2\times S^2$?
\end{question}

Research Institute for Mathematical Sciences, Kyoto University\\
Kitashirakawa Oiwake Sakyo-ku Kyoto-shi 606-8502\\
E-mail address: tange@kurims.kyoto-u.ac.jp

\end{document}